\documentclass[11pt]{amsart}
\usepackage[]{amsmath, amsthm, amsfonts, verbatim, amssymb,pb-diagram, stmaryrd, wasysym,tikz}
\usepackage{showkeys}
\usepackage{tikz}

\usepackage{epsfig}
\begin{document}
\newtheorem {theo}{Theorem}
\newtheorem {coro}{Corollary}
\newtheorem {lemm}{Lemma}
\newtheorem {rem}{Remark}
\newtheorem {defi}{Definition}
\newtheorem {nota}{Notation}
\newtheorem {prop}{Proposition}
\newtheorem {conj}{Conjecture}
\def\cD{{\mathcal D}}
\def\cE{{\mathcal E}}
\def\cM{{\mathcal M}}
\def\ocM{\overline{\mathcal M}}
\def\oM{\overline M}
\def\cF{\mathcal F}
\def\cH{\mathcal H}
\def\cO{\mathcal O}
\def\cG{\mathcal G}
\def\cI{{\mathcal I}}
\def\cR{{\mathcal R}}
\def\cG{{\mathcal G}}
\def\cl{{\mathcal L}}
\def\cL{{\mathcal L}}
\def\cT{{\mathcal T}}
\def\bR{{\mathbb R}}
\def\bZ{{\mathbb Z}}
\def\bQ{{\mathbb Q}}
\def\bC{{\mathbb C}}
\def\cm{{m}}
\def\ca{{a}}
\def\Alb{{\rm Alb}}
\def\cA{{\mathcal A}}
\def\cB{{\mathcal B}}
\def\cT{{\mathcal T}}
\def\cF{{\mathcal F}}
\def\cS{{\mathcal S}}
\def\cM{{\mathcal M}}
\def\cP{{\mathcal P}}
\def\cC{{\mathcal C}}
\def\cJ{\mathcal J}
\def\cX{\mathcal X}
\def\cU{\mathcal U}
\def\cY{\mathcal Y}
\def\cW{\mathcal W}

\def\spb{\smallpagebreak}
\def\mpb{\vskip 0.5truecm}
\def\bpb{\vskip 1truecm}
\def\wtM{\widetilde M}
\def\tM{\widetilde M}
\def\tN{\widetilde N}
\def\tC{\widetilde C}
\def\tX{\widetilde X}
\def\tY{\widetilde Y}
\def\bs{\bigskip}
\def\ms{\medskip}
\def\noi{\noindent}
\def\td{\nabla}
\def\pd{\partial}
\def\hol{$\text{hol}\,$}
\def\Log{\mbox{Log}}
\def\Log{\mbox{Log}}
\def\bP{\bf P}
\def\dxi{d x^i}
\def\dxj{d x^j}
\def\dyi{d y^i}
\def\dyj{d y^j}
\def\dzi{d z^I}
\def\dzj{d z^J}
\def\ozi{d{\overline z}^I}
\def\ozj{d{\overline z}^J}
\def\oz1{d{\overline z}^1}
\def\oz2{d{\overline z}^2}
\def\oz3{d{\overline z}^3}
\def\sI{\sqrt{-1}}
\def\hol{$\text{hol}\,$}
\def\ok{\overline k}
\def\ol{\overline l}
\def\oJ{\overline J}
\def\oT{\overline T}
\def\oS{\overline S}
\def\oV{\overline V}
\def\oW{\overline W}
\def\oI{\overline I}
\def\oK{\overline K}
\def\oL{\overline L}
\def\oj{\overline j}
\def\oi{\overline i}
\def\ok{\overline k}
\def\oz{\overline z}
\def\om{\overline mu}
\def\on{\overline nu}
\def\oa{\overline \alpha}
\def\ob{\overline \beta}
\def\of{\overline f}
\def\og{\overline \gamma}
\def\ogamma{\overline \gamma}
\def\odelta{\overline \delta}
\def\otheta{\overline \theta}
\def\ophi{\overline \phi}
\def\opd{\overline \partial}
\def\oA{\overline A} 
\def\oB{\overline B}
\def\oC{\overline C}
\def\op{\overline D}
\def\oIq1{\oI_1\cdots\oI_{q-1}}
\def\oIq2{\oI_1\cdots\oI_{q-2}}
\def\op{\overline \partial}
\def\ua{{\underline {a}}}
\def\us{{\underline {\sigma}}}
\def\Sym{{\mbox{Sym}}}
\def\bp{{\bf p}}
\def\bk{{\bf k}}
\def\obeta{\overline\beta}
\def\id{\mbox{id}}
\def\tq{\widetilde{q}}
\def\ti{\widetilde{i}}
\def\tGamma{\widetilde{\Gamma}}
\def\tV{\widetilde{V}}
\def\CM{\mbox{CM}}
\def\rank{\mbox{rank}}
\def\ocA{\overline{\cA}}
\def\hA{\widehat{\cA}}
\def\hM{\widehat{M}}
\def\tr{\mbox{Tr}}
\def\ad{\mbox{ad}}
\def\nni{\noindent}
\def\oZ{\overline Z}
\def\te{\widetilde e}
\def\IIm{\mbox{Im}}
\def\oY{\overline Y}
\def\Exp{\mbox{Exp}}
\def\ad{\mbox{ad}}
\def\Rre{\mbox{Re}}
\newcommand{\hooklongrightarrow}{\lhook\joinrel\longrightarrow}

\noi
\title[Open Torelli locus and complex ball quotients]
{Open Torelli locus and complex ball quotients}
\author[Sai-Kee Yeung]
{Sai-Kee Yeung}

\begin{abstract}  
{\it  We study the problem of non-existence of totally geodesic complex ball quotients in the open Torelli locus
in a moduli space of principally polarized Abelian varieties using analytic techniques. }
\end{abstract}

\address[]{Mathematics Department, Purdue University, West Lafayette, IN  47907, USA} \email{yeung@math.purdue.edu
}

\thanks{\noindent{
The author was partially supported by a grant from the National Science Foundation}}

\noi{\it }
\maketitle

\bigskip

\begin{center}
{\bf \S1. Introduction}
\end{center}

\ms
\noi{\bf 1.1}  Let $\cM_g$ be the moduli space or stack of Riemann surfaces of genus $g\geqslant 2$.  Let $\ocM_g$ be the Deligne-Mumford compactification of
$\cM_g$.  Let $\cA_g$ be the moduli space of
principally polarized Abelian varieties of complex dimension $g$.   We know that $\cA_g=\cS_g/Sp(2g,\bZ)$ is the quotient of
 the Siegel Upper Half Space $\cS_g$ of genus $g$.  Let $\ocA_g$ be the Bailey-Borel compactification of $\cA_g$.  Associating a smooth Riemann surface
 represented by a point in $\cM_g$ to its Jacobian, we obtain the Torelli map $j_g:\cM_g\rightarrow \cA_g$.  The Torelli map extends to
 $j_g:\ocM_g\rightarrow \ocA_g$.
 The image $T_g^o:=j_g(\cM_g)$ is called the open Torelli locus of $\cA_g$.    It is well-known that the Torelli map $j_g$ is injective on $\cM_g$.  As a mapping between stacks, the mapping $t_g|_{\cM_g}$ is known
 to be an immersion apart from the hyperelliptic locus, which is denoted by 
$H_g$.   $H_g$ is the set of points in $\cM_g$ parametrizing hyperelliptic  curves of genus $g$.

It is a natural problem to study $\cM_g$, $j_g$ and to characterize the Torelli locus in $\cA_g$.  There are many interesting directions and approaches
to the problems.  
Our motivation comes from the following conjecture in the literature.  

\begin{conj} (Oort [O])
Let $T_g^o$ be the open Torelli locus in the Siegel modular variety $\cA_g.$   Then for $g$ sufficiently large,
the intersection of $T_g^o$  with any Shimura variety $M\subset \cA_g$ of strictly positive dimension is not
Zariski dense in $M$.
\end{conj}

The problem is related to a conjecture of Coleman [C] that the cardinality of CM points on $\cM_g$ cannot be infinite if $g$ is sufficiently large.
Shimura varieties are arithmetic locally Hermitian symmetric spaces.   Hence we may consider a geometrically slightly more general question of
whether there exists a locally Hermitian symmetric space in $\cA_g$ with a Zariski open set in $T_g^o$.  In such case, the lattice $\Gamma$
involved in the 
complex rank one case, namely complex balls $B^n_{\bC}=PU(n,1)/P(U(n)\times U(1))$, may not be arithmetic.  

Conjecture 1 is open, but there are quite a few interesting partial results.  First of all there is the result of Hain [Ha] that for $M$ locally Hermitian symmetric of rank at least $2$, it cannot happen that $M\subset T_g(\ocM_g)-(t_g(H_g)\cup t_g(\ocM_g\backslash \cM_g))$, cf. also de Jong-Zhang [dJZ] for 
precise formulation and some results.  On the other hand,
there is the result of de Jong-Noot [dJN] that there are examples of Shimura curves in $\cM_g$ for $g=4$ and $6$.  A more systematic and
complete treatment in this direction can be found in Moonen [Moo].
 To the knowledge of the author,
not much is known for complex ball quotients of
dimension $n\geqslant 2$, apart from some restrictions in terms of Higgs bundles given in Chen-Lu-Tan-Zuo [CLTZ].  A possible reason is that on the one
hand rigidity results in general are not strong enough for super-rigidity properties in the rank one complex cases, and  on the other hand the problem
is not concrete enough to be handled by geometric techniques developed for specific Riemann surfaces.

\ms
\noi{\bf 1.2}
The goal of this paper is
to provide a method  which is applicable to locally Hermitian symmetric spaces and in particular to all
complex ball quotients of
dimension $n\geqslant 2$.   It is a weaker statement but is in support of Conjecture 1. 

\begin{theo}  
 The set $T_g^o- j_g(H_g)\subset \cA_g$ for $g\geqslant 2$ does not contain any  complex hyperbolic complex ball quotient, compact or non-compact with finite volume,
  of complex dimension
at least $2$ as a totally geodesic complex suborbifold of $\cA_g$.
\end{theo}

\ms
\noindent {\bf Remarks}  \\
The method of proof applies immediately to other locally Hermitian symmetric spaces of complex dimension at least $2$.  We refer the readers to Section {\bf 3.5} and {\bf 4.4} for more details.  However, the
end result for $\rank_{\bR}M\geqslant 2$ followed already from the results of [Ha], where the proof is completely different.

\bs
Combining Theorem 1, the results of [Ha], the results of [M\"o] and the very recent result of [AN], we have now a rather complete 
picture for Shimura varieties on the complement of the hyperelliptic locus in the open Torelli locus.

\begin{theo}  
Let $g\geqslant 2$.
The space $T_g^o- j_g(H_g)\subset \cA_g$ in the Siegel modular variety $\cA_g$ does not contain any Shimura subvariety of $\cA_g$,
except when $M$  is the Torelli image of a Riemann surface with genus $g=3,4$.  There is
 only one curve for each of $g=3, 4$, with universal families given by $y^4=x(x-1)(x-t)$ and $y^6=x(x-1)(x-t)$ respectively.
\end{theo}

It is a pleasure for the author to thank Ngaiming Mok for very helpful comments and suggestions.
The author is also indebted to the referee for pointing out mistakes in earlier drafts of the paper and for other useful comments.
 
\bs

\begin{center}
{\bf \S2. Preliminaries and rigidity}
\end{center}

\noi{\bf 2.1}  The approach we take is complex analytic, trying to compare various Kobayashi metrics making use of Schwarz Lemma and results in rigidity.
The reader may refer to {\bf 3.1} for a brief summary of facts needed about Kobayashi metric.

As mentioned in the introduction, we may consider either $M$ as a smooth submanifold of $\cA_g$, or a suborbifold.  In the latter case the Kobayashi metrics are considered to
be in orbifold sense as follows.  Recall that all the singularities of $\cA_g$ are orbifold singularities since
$\cA_g$ is a quotient of the Siegel Upper Half Space $\cS_g$ by a discrete group and $\cS_g$ is smooth.   By an orbifold embedding (resp. mapping) $f:M\rightarrow \cA_g$, we mean that there is
a finite covering $\pi:\cA_g'\rightarrow \cA_g$ so that $\cA_g'$ is smooth and there is an embedding (resp. mapping) $f':M'\rightarrow \cA'$ for which $M'$ is smooth and $f\circ \pi=\pi\circ f'$.
Since Kobayashi metric is invariant under a biholomorphism and in particular invariant under a local holomorphic covering map, the argument throughout the article would be 
independent of the local uniformization taken at each orbifold singularity.

Hence the Kobayashi metric studied throughout the article is in orbifold sense as explained.

\ms
\noi{\bf 2.2}  The  proof of Theorem 1
makes use of  the following result in [A].  

\begin{prop} 
 ([A]  Theorem 1.1)
Let $\tM=B_{\bC}^N$ be the complex unit ball of dimension $N\geqslant 2$.  There is no holomorphic embedding of  $\tM$ into $\cT_{g}$ which is isometric with 
respect to the Kobayashi metrics on $\tM$ and $\cT_g$.
\end{prop}

\bs

\begin{center}
{\bf \S3. Totally geodesic embeddings from complex balls to Siegel Upper spaces}
\end{center}

\bs
\noi{\bf 3.1} Our approach relies on basic properties concerning totally geodesic embeddings of complex balls in the Siegel Upper Half Space.  The purpose of this section
is to explain results in this direction relevant to our purpose.  For this purpose, it is more
convenient to consider the bounded model of $\cS_g$, namely, the classical bound domain $III_g$.  The main technical result of this paper is Proposition 3 stated in {\bf 3.4}.

Basic properties of classical domains can be found in [He], [Mok1] and [Sat].  Since we are going to use the results of Satake, we follow closely 
the exposition in [Sat].  First let us explain briefly classical domains of type $I_{p,q}$ and $III_k$ according to {\bf 1.2, 1.3} of [Sat], of which the terminology is to be used
in later parts of this section. 

\bs

\noi{\it (i). $I_{p,q}$:}  Consider $V$ a vector space over $\bC$ equipped with a non-degenerate Hermitian form $F$ of signature $(p,q)$ with $p>q>0$.
$I_{p,q}=\cD(V,F)=:\cD$ is the space of $q$ dimensional complex subspace $V_-$ of $V$ so that $F|_{V_-}$ is negative definite.  Let 
$V_+$ be the orthogonal complement of $V_-$ in $V$, so that $F|_{V_+}$ is positive definite.  A point in $\cD$ is determined
by $V_-$, or the pair $(V_+, V_-)$.  Let  $z_0$ be a fixed point $\cD$, determined by $(V_+^{(0)},V_-^{(0)})$ and for convenience can be chosen to be the origin.  Let $(e_1,\cdots,e_p)$ and
$(e_{p+1},\cdots,e_{p+q})$ be orthonormal basis of $V_+^{(0)}$ and $V_-^{(0)}$ respectively, so that together they form a basis
of $V$.  A point $z\in \cD$ is now determined by $(V_+, V_-)$ with $V_-$ spanned by the basis 
\begin{equation}
\sum_{i=1}^pe_iz_{ij}+e_{p+j}, \ 1\leqslant j\leqslant q, 
\end{equation}
where the $(p,q)$ matrix $Z=(z_{ij})$ satisfies
$I_q-^t\oZ Z>0$.   Denote by $M_{p,q}$ the space of all $p \times q$ matrices with entries in $\bC$.   Identifying $\cD$ with $\{Z\in M_{p,q}:I_q-\oZ^t Z>0\}$, we realize $I_{p,q}$ as a bounded domain in $\bC^{pq}$.

The complex ball $B_{\bC}^n$ in $\bC^n$ is just $I_{n,1}$.  

\bs

\noi{\it (ii).  $III_k$:}  Consider $V_{\bR}$ a vector space over $\bR$ of dimension $2k$ equipped with a non-degenerate alternating bilinear form $A$.
Let $V=V_{\bC}$ be the complexification of $V_{\bR}$.  $A$ extends naturally to $V$.  The Hermitian form defined by
\begin{equation}
F(x,y):=iA(\bar x,y)
\end{equation}
has signature
$(k,k)$ on $V$.  $III_k=\cD(V_{\bR},A)$, or simply, $\cD$, is the space of all complex structures $I$ on $V_{\bR}$ so that the bilinear form $A(x,Iy)$ is symmetric
and positive definite.  Let $W=\{x\in V|Ix=ix\}$ so that $V=W+\overline W$.  It follows that 
\begin{equation}
A|_W=0, \ \ F|_W>0
\end{equation}
and $\overline W$ is the orthogonal
complement of $W$ in $V$ with respect to $F$.  Hence $I\in \cD$ is determined by $W$ or the pair $(W,\overline W)$ satisfying (3).  Fix a point
$z_o\in \cD$ corresponding to $W^o$.  Let $(e_1,\cdots,e_p)$ be an orthonormal basis of $W^o$ with respect to $F$ and let 
\begin{equation}
e_{k+i}=\overline e_i, \ 1\leqslant i\leqslant k.
\end{equation}  As described in $I_{p,q}$ with $p=q=k$ in the description of $W\in \cD$, it follows that $W$ is described by a
$k\times k$ symmetric 
complex matrices
$Z$ with $I_k-\oZ^t Z>0$.  Hence
$\cD$ is identified with the bounded domain $\{Z\in M_{k,k}:Z=Z^t, I_k-\oZ^t Z>0\}\subset \bC^{n(n+1)/2}$.  

The Siegel Upper Half Space
$\cS_g$ is biholomorphic to $III_g$.

\bs
\noi{\bf 3.2} Here we recall briefly the classification of holomorphic totally geodesic embeddings of a Hermitian symmetric domain into another.
In general, the  classification of holomorphic totally geodesic embedding of a Hermitian symmetric space $N_1=G_1/K_1$ into another Hermitian symmetric space
$N_2=G_2/K_2$ with respect to the 
Bergman metrics has been given by Satake [Sat] and Ihara [I], where $G_i$ is a semi-simple Lie group and $K_i$ a maximal compact subgroup for $i=1,2$.
 Since the manifolds involved are symmetric, the classification of totally geodesic embeddings is reduced to the classification of injective Lie algebra homomorphisms
 $\rho:\mathfrak{g_1}\rightarrow \mathfrak{g_2}$ for the corresponding Lie groups.
 The invariant complex structure on $N_i$ is given by an element $H_{oi}\in K_i$.  The totally geodesic embedding is holomorphic if the condition $(H_1)$, namely,
 $\rho\circ \ad(H_{o1})=\ad(H_{o2})\circ \rho$ is satisfied.  The condition $(H_1)$ is a consequence of the condition $(H_2)$, namely, $\rho(H_{o1})=H_{o2}$.
 This is explained on page 427 of [Sat].  
   Some explanation in terms of root system and Dynkin diagrams  by Ihara in [I].   We refer the reader to [Sat] for any unexplained notation and terminology.

\bs
\noi{\bf 3.3} We consider now the specific situation of classification of holomorphic totally geodesic embeddings of $I_{N,1}$ into $III_g$.  Let $0$ be the origin in $B^N_{\bC}$, which may be assumed to be mapped to the origin $0$ of $\cS_g$ realized as the bounded domain $III_g$
as above, since the spaces involved are homogeneous.
The classification is described in Theorem 1, {\bf 3.2} and the Table on page 460 of [Sat].  From the description there, any holomorphic embedding of $I_{N,1}$ to $III_g$ is 
a  direct sum of a number of compositions of
the following types of totally geodesic mappings.  

\ms  
\noi Type 1,  standard embeddings: This includes embedding $i_{1a}: I_{N,1}\rightarrow I_{p,q}$, for $N\leqslant p$ and $i_{1b}: III_{k}\rightarrow III_{l}$, for $k<l$,  given by the standard representation or standard embedding. 

\ms
\noi Type 2, connecting embeddings:  The embedding  $i_2=\iota_{p,q}: I_{p,q}\rightarrow III_{p+q}$ is given in page 432-433 of [Sat],
$$I_{p,q}\ni Z\longmapsto \left(\begin{array}{cc}
  0&^tZ\\
  Z&0
  \end{array}
  \right)\in III_g\cong \cS_{p+q}.$$

\ms
\noi Type 3, absolutely irreducible embeddings:   The embedding $i_3:I_{p,1}\rightarrow I_{r,s}, r={p\choose m}, s={p\choose m-1}$, and $i_3 :I_{p,1}\rightarrow I_r$ in the case of $r=s$, which happens when
 $p\equiv 1\pmod 4$ and $m=\frac{p+1}2$.
 This corresponds to the skew-symmetric tensor representations of degree $m$ as explained in {\bf 3.2} and the Table on page 460 of [Sat].    
 
The construction is described in page 448 of [Sat].  The representation is given by $\rho=\Lambda_m$ of $G\cong SU(p,1)$, corresponding to skew-symmetric tensors of degree $m$.  In 
terms of a fundamental system of roots of the Lie algebra involved, the highest weight $\lambda_\rho$ of $\rho$ is given by
$$\lambda_\rho=(\overbrace{1,1,\cdots,1}^\text{m},0,\cdots 0)$$
with $1\leqslant m\leqslant p$.  Recall the setting given in  {\bf (3.1)(i)} in the case of $I_{p,1}$.
Let $(e_i)_{i=1,\cdots,p+1}$ be an orthonormal basis of $V$, a basis of the exterior algebra $\Lambda^m$ is given by $e_{i_{1}\cdots  i_m}=e_{i_1}\wedge\cdots\wedge e_{i_m}$, $i_1<i_2<\cdots<i_m$.
As in {\bf (3.1)(i)}, there is a Hermitian form $F$ on $V$, which induces $F^{(m)}$ on $\Lambda^m$ given by 
$$F^{(m)} (x_1\wedge\cdots\wedge x_m,y_1\wedge\cdots\wedge y_m)=\det (F(x_i,y_j))$$
for $x_i, y_j\in V$ and is invariant under $\rho(\mathfrak{g}_{su_{p,1}})$.  Since
$$F^{(m)}(e_{i_1\cdots  i_m},e_{i_1\cdots  i_m})=\left\{\begin{array}{ll}
1& \ i_m<p+1, \\
-1& \ i_m= p+1,
\end{array}\right.
$$
the Hermitian for  $F^{(m)}$ has signature  $(r,s)$ with $r={p \choose m}, s={p \choose m-1}$.  

The totally geodesic isometry of the symmetric domains $f:I_{p,1}\rightarrow I_{r,s}$ 
described here is given in page 448 of [Sat] by 

\begin{equation}
\begin{array}{ccccccc}
\cD(V,F)&\ni& (V_+,V_-)&\longmapsto&(\Lambda_m(V_+)\otimes 1, \Lambda_{m-1}(V_+)\otimes V_-)&\in& \cD(\Lambda_m(V),F^{(m)}).\\
\parallel&&\parallel&&\parallel&&\parallel\\
I_{p,1}&\ni&z&\longmapsto&z'&
 \in &I_{r,s}
\end{array}
\end{equation}

\ms
In the case of representation in $III_r\subset I_{r,r}$, this corresponds to the above discussion with 
$r=s$ and hence $m=\frac{p+1}2$ and  $p\equiv 1\pmod 2$.   In such case, we can define a Bilinear form $B$ on $\Lambda_m(V)\times\Lambda_m(V)$ by
$$x\wedge y=B(x,y)e_{1\cdots p+1}$$
which satisfies 
$$B(y,x)=(-1)^{m^2}B(x,y).$$
Hence if $m\equiv 1\pmod 2$ or 
$p\equiv 1\pmod 4$, the bilinear form $B(x,y)$ is  an alternating bilinear form.  Furthermore, 
there is a semi-linear transformation $\sigma$ on $\Lambda_m(V)$ so that
$$F^{(m)}(x,y)=iB(x^\sigma,y),$$
where $\sigma$ satisfying $\sigma^2=1$ is explicitly written in [Sat], page 449,  as follows.  Let $M=(i_1,\cdots,i_m)$ be an oriented subset of $(1,2,\cdots,p+1)$ and
$M^c$ the complement.  Then
\begin{equation}
e_M^\sigma=a(M)e_{M^c}, \ a(M)=-i\epsilon(M^c,M)\eta(M),
\end{equation}
where $\epsilon(M^c,M)=\pm 1$ is the signature of the permutation of $(M^c, M)$ with respect to $(1,2,\cdots,p+1)$ and $\eta(M)=-1$ (resp. $1$) if $p+1\in M$ (resp. $p+1\not\in M$).
Hence $a(M)=\pm i$.  $\sigma$ serves and complex conjugate as in {\bf 3.1(ii)}.
In this case, the bilinear form $B$ serves as the bilinear form $A$ as needed in equation (2) in
{\bf 3.1(ii)} for the definition of $III_r$.
The totally geodesic isometry of the symmetric domains $f:I_{p,1}\rightarrow III_r$  is given by (5) with $r=s$.

 We summarize the result of [Sat], which is relevant to us, from Theorem 1 and the Table on page 460 in [Sat].  As explained in {\bf 3.2}, totally geodesic
holomorphic embedding corresponds to condition $(H_1)$, which by (a) below reduces the problem to representations satisfying condition $(H_2)$.  The Table on page 460 in [Sat]
summarizes the representations satisfying condition $(H_2)$ determined in Section 3 of  [Sat].

\begin{prop}(Satake)\\
(a). Let $\rho$ be a representation of $\mathfrak{g}=\mathfrak{su(m,1)}$ into $III_g$ satisfying condition $(H_1)$.  Then there exists
absolutely irreducible representations $\rho_i$ $(1\leqslant i\leqslant r_1)$ of $\mathfrak{g}$ into $(III)_{p_i}$ $(p_i>0)$ satisfying $(H_2)$ and 
absolutely irreducible representations $\rho_i$ $(r_1\leqslant i\leqslant r_1+r_2)$ of $\mathfrak{g}$ into $(I)_{p_i,q_i}$ 
$(p_i, q_i\geqslant 0, p_i+q_i>0)$ satisfying $(H_2)$ such that 
$\rho$ is $k$ equivalent to the direct sum of representations $\sum_{i=1}^{r_1}\rho_i+\sum_{i=r_1+1}^{r_1+r_2}\iota_{p_i,q_i}\circ\rho_i$ up to a trivial representation, where
$\sum_{i=1}^{r_1}p_i+\sum_{i=r_1+1}^{r_2}(p_i+q_i)\leqslant g$.\\
(b). An absolutely irreducible representation $\rho$ of $\mathfrak{g}$ into $I_{r,s}$ or $III_r$ satisfying $(H_2)$ corresponds to embeddings of type $i_3$ described earlier.
\end{prop}

We refer the reader to the original source [Sat] for any unexplained terminology.    In particular, the notion of $k$-equivariant  and direct sum are
described in \S1 and \S2 of [Sat].  The result of Satake applies to
any semi-simple Lie algebra $\mathfrak{g}$ of Hermitian type.



\ms
Let us now describe the mapping given in (5) more carefully, which is to be used later.  

\bs

\begin{lemm}
In terms of standard coordinates, the totally geodesic mapping $f:B_n\rightarrow I_{r,s}\subset\bC^N$ or $III_r\subset\bC^N$ for $r=s$ with $f(0)=0$ with
$z'=f(z)$ as described above  is linear in $z$, with image given by the intersection $f(B_n)$ with a subspace of $\bC^N$ of appropriate dimension.  
\end{lemm} 

\ms
\noi{\bf Proof}   We remark that the standard coordinates as used in [Sat] in the description above are also the Harish-Chandra coordinates.

Consider first $f:I_{p,1}\rightarrow I_{r,s}$ as given by (5).
 In terms of (1) for $I_{p,1}$, the point $z$ in (5) corresponds to $V_+$ being spanned by 
$\sum_{i=1}^p e_iz_i+e_{p+1}$.  Since it is a holomorphic totally geodesic embedding, the mapping is equivariant with respect to the action of
$G$ and in particular invariant under the action of the isotropy group $K=S(U(p)\times U(1))$.  In particular, it suffices for us to investigate the image 
$f(z)$ for $z=(z_1,0,\cdots,0)$ with $|z_1|^2<1$.   Again, we use $(e_1,\cdots, e_p)$ and $(e_{p+1})$ to denote some orthonormal basis of $V_+^{(0)}$ and $V_-^{(0)}$ respectively.
With $z$ as described, an orthonormal basis of $V_-$ and $V_+$ at $z$ are  $e_{p+1}'$ and $(e_1',e_2,\cdots,e_p)$ respectively, where 
\begin{equation}
e_{p+1}'=\frac1{\sqrt{1+|z|^2}}(z_1e_1+e_{p+1}), \ 
e_1'=\frac1{\sqrt{1+|z|^2}}(e_1-\oz_1e_{p+1}).
\end{equation}

To describe $z'$ in the image of $f$ in (5), we need to investigate $\Lambda_{m-1}(V_+)\otimes V_-$  where $F^{(m)}$ is negative definite,
and express them in terms of base vectors of $\Lambda_m^{(0)}=\Lambda_m(V_+^{(0)})\oplus \Lambda_{m-1}(V_+^{(0)})\otimes V_-^{(0)}$.
From definition, $\Lambda_{m-1}(V_+^{(0)})\otimes V_-^{(0)}$ is generated by a basis consists of the following two types of
elements,\\
(i) $e_{i_1\cdots i_{m-1}}\wedge e_{p+1}'$, where $1<i_1<\cdots<i_{m-1}\leqslant p$,\\
(ii) $e_1'\wedge e_{i_1\cdots i_{m-2}}\wedge e_{p+1}'$,  where $1<i_1<\cdots<i_{m-2}\leqslant p.$\\

From the formula of $ e_{p+1}'$ in (7), we compute in case (i) that
\begin{equation}
e_{i_1\cdots i_{m-1}}\wedge e_{p+1}'=\frac1{\sqrt{1+|z|^2}}((-1)^{m-1}z_1e_{1i_1\cdots i_{m-1}}+e_{i_1\cdots i_{m-1}(p+1)})
\end{equation}
which is proportional to 
\begin{equation}
(-1)^{m-1}z_1e_{1i_1\cdots i_{m-1}}+e_{i_1\cdots i_{m-1}(p+1)}.
\end{equation} 
Notice that we need the coefficient of $e_{i_1\cdots i_{m-1}(p+1)}$ to be $1$ in the format of (1).

In case (ii), (7) gives 
$$e_1'\wedge e_{p+1}'=e_1\wedge e_{p+1}$$ 
and hence
\begin{equation}
e_1'\wedge e_{i_1\cdots i_{m-2}}\wedge e_{p+1}'=e_{1i_1\cdots i_{m-2}(p+1)}.
\end{equation}

 Let $u={p-1\choose m-1}$.
It follows from the above explicit computation that for $f:I_{p,1}\rightarrow I_{r,s}$ as given in (5), the coordinates of $z'=f(z)$ are given by 
\begin{equation}
z=(z_1,0,\cdots,0)\mapsto f(z)=z'  \  \mbox{with}  \ z'_{ij}=\left\{\begin{array}{cc}
(-1)^{m-1}z_1,& 1\leqslant i\leqslant t, 1\leqslant j\leqslant u\\
0& \mbox{otherwise}
\end{array}
\right.
\end{equation}

Similar constructions apply to  $f:I_{p,1}\rightarrow III_r\subset I_{r,r}$ corresponding  (5) with $r=s$.
We recall that the coordinates $z'$ in $III_r$ in the image of $f$ is determined according to (1) in {(3.1)} with respect to a corresponding basis of vectors in 
$(\Lambda_m(V_+)\otimes 1, \Lambda_{m-1}(V_+)\otimes V_-)$.  Moreover, the choice of the base vectors of is given by (4),  choosing $e_{r+i}$ to
be the complex conjugate of $e_i$.
In our case of $\Lambda_m^{(0)}$, the complex conjugate 
 is given by $\sigma$ in the setting of (5) and $e_{I_m}^\sigma=a(I_m)e_{I_m^c}$ in terms of earlier notations,
here $I_m=i_1\cdots i_m$ is an index set.
Hence $z'=(z'_{ij})_{1\leqslant i,j\leqslant r}$ are determined by having a basis of $\Lambda_{m-1}(V_+)\otimes V_-$ of the form
\begin{equation}
\sum_{i=1}^r e^\sigma_{I_i(p+1)}z'_{ij}+e_{I_j(p+1)}=\sum_{i=1}^r a(I_i(p+1))e_{(I_i(p+1))^c}z'_{ij}+e_{I_j(p+1)},\  1\leqslant j\leqslant r
\end{equation}
where $a(I_i(p+1))=\pm1$.  From the expressions in (9) and (10), we see that (11) still applies in the sense that $z'_{ij}$ are either $0$ or $(-1)^m z_1$.    This computation shows that in terms of the standard coordinates of the bounded domains in $\bC^N$ as described in {\bf 1.3, 1.4} of [Sat] or Chapter 4 \S2 of [Mok 1] for the classical domains
that the image is the intersection of a line in $\bC^n$ with $I_{r,s}$ or $I_r$.
As mentioned earlier, this works for any complex direction obtained under the action of the isotropy group at $0$ on the domain.  Since the mapping $f$ is a totally
geodesic mapping and is hence equivariant under the action of the isotropy group,  the lemma follows.

\qed

\bs
\noi{\bf 3.4}  We now state the main result of this section.

 \begin{prop} 
Let  $i:B^N_{\bC}\rightarrow III_g$ be a totally geodesic embedding.  Then there exists a holomorphic map $p:III_g\rightarrow B^N_{\bC}$ so that $p\circ i$ is the identity mapping on $B^N_{\bC}$.
\end{prop}

\noi{\bf Proof}  Again, we let $0$ be the origin in $B^N_{\bC}$, which may be assumed to be mapped to the origin $0$ of $III_g$, a bounded domain realization of
$\cS_g$,
as discussed earlier.   To streamline the presentation, let us consider each simple type of presentations $i_1, i_2, i_3$ in details before the general case described in Proposition 2.


\ms
\noi  Type 1:  This includes  embedding is given by $i_1: I_{N,1}\rightarrow I_{p,q}$, for $N\leqslant p$ or $i_{1b}: III_{k}\rightarrow III_{l}$, for $k<l$, given by the standard embedding into the corresponding  upper left hand corner of the image.

The classical domain  $I_{p,q}$ is given as a symmetric space $G/K$ with $G=SU(p,q)$ and $K=S(U(p)\times U(q))$, where $K$ is the isotopy group at $0$.
For a holomorphic totally geodesic
embedding  $i_{1a}:B^N_{\bC}\cong I_{N,1}\rightarrow I_{p,q}$, conjugating by $K$ if necessary, we may assume that $\frac{\pd}{\pd z_{11}}\in (i_1)_*(T_{B^N_{\bC}})$.  Now we observe that
$\frac{\pd}{\pd z_{ij}}$ with $i>1$ and $j>1$ cannot lie in $(i_{1a})_*(T_{B^N_{\bC}})$, for otherwise the image $i_{1a}(B^N_{\bC})$ as a symmetric space would have real rank at 
least $2$ considering the tangent vectors $\frac{\pd}{\pd z_{11}}$ and $\frac{\pd}{\pd z_{ij}}$, contradicting the fact that $I_{1,N}$ has real rank $1$.
  It follows that $i_{1a}(B^N_{\bC})$ has to lie in one of the following two subspaces of $I_{p,q}$
$$I_{1,q}=\{z=[z_{ij}]\in I_{p,q}|z_{ij}=0 \ \mbox{for} \ i\geqslant 2\} \ \mbox{or}\ I_{p,1}=\{z=[z_{ij}]\in I_{p,q}|z_{ij}=0 \ \mbox{for} \ j\geqslant 2\}.$$ In either case, $N\leqslant \max(p,q)=q.$

Suppose $i_{1a}(B^N_{\bC})\subset I_{1,q}$.   Conjugating by some elements in $K$ if necessary, we may assume that mapping $i_{1a}$ is given by 
$$i_{1a}(z_1,\cdots,z_N)=\left[\begin{array}{cccccc}
z_1&\cdots&z_N&0&\cdots&0\\
0&\cdots&0&0&\cdots&0\\
\vdots&&&&&\vdots\\
0&\cdots&0&0&\cdots&0
\end{array}\right].
$$
Consider the holomorphic projection map $p_1:I_{p,q}\rightarrow B_{\bC}^N$ given by
$$p_1(\left[\begin{array}{cccccc}
z_{11}&\cdots&z_{1N}&z_{N+1}&\cdots&z_{1q}\\
\vdots&&&&&\vdots\\
z_{p1}&\cdots&z_{pN}&z_{N+1}&\cdots&z_{pq}
\end{array}\right])
=\left[\begin{array}{cccccc}
z_{11}&\cdots&z_{1N}&0&\cdots&0\\
0&\cdots&0&0&\cdots&0\\
\vdots&&&&&\vdots\\
0&\cdots&0&0&\cdots&0
\end{array}\right].
$$
Denote by $Y$ the matrix in the domain and $Z$ the matrix in the image.
As the $p\times q$ matrix $Y$ satisfies $I-^t\oY Y>0$ from definition, it follows that
$\sum_{i=1}^N|z_i|^2\leqslant \sum_{i=1}^p|z_i|^2<1$ and hence the image of $p_1$ lies in $B_{\bC}^N$.
Furthermore, it follows from definition that $p_1\circ i_{1a}=1_{B_{\bC}^N}.$  Hence $p_1$ gives us the retraction that we need.

 For $i_{1b}: III_{k}\rightarrow III_{l}$, for $k<l$, given by the standard embedding into the corresponding  upper left hand corner of the image.  Denote by $p_{1b}:III_{l}\rightarrow III_{k}$ 
  the projection onto the upper left hand corner
 \begin{eqnarray*}
p_{1b}([z_{ij}]_{i,j=1,\cdots,l})&=&[z_{ij}]_{i,j=1,\cdots,k}.\\
\end{eqnarray*}
Let $Y=[z_{ij}]_{i,j=1,\cdots,l}$ and $U=[z_{ij}]_{i,j=1,\cdots,k}$.  The fact that $Y$ is symmetric implies that $U$ is symmetric.
Now $I_g-\overline{Y}Y>0$ implies that $I_{g'}-\overline{U}{U}>0$ as each column vector of $U$ is part of a column vector of $Y$.  Hence the image lies in
$III_k$.

It is clear that $p_{1b}\circ i_{1b}|_{III_{g'}}$ is the identity map and hence $p_3$ is a projection.

\ms
\noi  Type 2:  $i_2=\iota_{p,q}:I_{p,q}\rightarrow III_{p+q}$ with $\iota_{p,q}(Z)=\left(\begin{array}{cc}
  0&^tZ\\
  Z&0
  \end{array}
  \right).$
Define 
$i_2:III_{p+q}\rightarrow I_{p,q}$ the projection

$${III_{p+q}}\ni Y=\left(\begin{array}{cc}
  W_1&^tZ\\
  Z&W_2
  \end{array}
  \right)\longmapsto \left(\begin{array}{cc}
  0&^tZ\\
  Z&0
  \end{array}
  \right)\longmapsto  Z\in {I_{p,q}}.$$
 Here $Y$ is symmetric.  From the fact that $I_{p+q}-\overline{Y}Y>0$, it follows that 
 $I_q-^t\overline{W} W-^t\overline{Z}Z>0$ and hence $I_q-^t\overline{Z}Z>0$.  Hence the image of 
 $i_2$ is really in $I_{p,q}$.
It follows from definition that $p_2\circ i_2=1_{I_{p,q}}$ is the identity map.

\ms
\noi  Type 3:  $i_3: I_{p,1}\rightarrow I_{r,s}$ or $III_r$ corresponding to the skew-symmetric tensor representations of degree $m$.  Since the second case can be considered
as a special case of  the first case, it suffices for us to consider the case of $I_{r,s}$ being the image.
From Lemma 1 in {\bf 3.3}, the image of $i_3$ is
a linear subspace $R:=\IIm(i_3)$ passing through the origin in $III_r$.  We claim that the argument of Lemma 1 shows that 
there is a projection $p_3:I_{r,s}\rightarrow f(I_{p,1})$ so that $p_3\circ i_3|_{I_{p,1}}$ is identity on $I_{p,1}$.  We actually take $p_3$ be the orthogonal projection of $I_{r,s}$ to the complex linear subspace $P$ of $\bC^N$ 
containing $f(I_{p,1})$ as a subdomain, after Lemma 1
To prove the claim, since the domain involved is a convex domain in $\bC^N$, 
it suffices for us to show that $p_3(I_{r,s})=f(I_{p,1})$.   Clearly $p_3(I_{r,s})$ contains $f(I_{p,1})$.
This in turn follows if we can prove  the corresponding statement for the projection of $I_{r,s}$ to $f(B^1_{\bC})$ for a geodesic 
$B^1_{\bC}\subset B^p_{\bC}\cong I_{p,1}$ through the origin, since $f$ is equivariant with respect to the action of the isotropy groups at $0$.
Hence it suffices for us to show that the projection $p_4$ of $I_{r,s}$ to the line in $\bC^N$ containing $f(B_{\bC}^1)$ is actually $f(B_{\bC}^1)$, where  
$B^1_{\bC}= \{(z_1,0,\dots,0):|z_1|<1\}$ as studied in the proof of Lemma 1.

Recall that $r={p\choose m}, s={p \choose m-1}$ and $u={p-1\choose m-1}$.  From (11) in the proof of Lemma 1, the image of $f$  is \begin{equation}
f(B^1_{\bC})=\Big\{\left(\begin{array}{ll}
z_1I_u & 0_{u,a}\\
0_{b,u}&0_{b,a}
\end{array}\right), \ (z_1,0,\dots,0)\in B^1_{\bC}\Big\}
\end{equation}
where $b={p-1\choose m}, a={p-1\choose m-2}$, $I_t$ is the identity matrix of size $t$ and $0_{c,d}$ is the zero matrix of size $c\times d$.  Clearly 
$$f(B^1_{\bC})\subset I_{r,r}\cong \left(\begin{array}{ll}
I_{r,r}& 0_{u,a}\\
0_{b,u}&0_{b,a}
\end{array}\right),
$$
where $I_{r.r}$ is a bounded symmetric domain of type I.

Similar to the construction of $I_{1a}$ in the standard embedding earlier, it is clear that there is a projection $q_1:I_{r,s}\rightarrow I_{u,u}$ as explained 
above by taking zeros in non-relevant entries.  Hence it suffices for us to show that there is a retraction $q_2: I_{u,u}\rightarrow \{(z_1I_u):|z_1|<1\}$.  Let $w=(w_{ij})\in I_{u,u}$.  By definition, it
satisfies
$I_q-\bar w^tw>0$ and hence $|w_{ij}|<1$ for each $i,j$. 
It suffices for us to define 
$$q_2(w)=(\frac1{u^2}\sum_{i,j=1}^u w_{ij})I_u.$$  Clearly $|\frac1{u^2}\sum_{i,j=1}^u w_{ij}|<1$ as $|w_{ij}|<1$.  Furthermore $q_2(\zeta I_u)=\zeta I_u$ for $\zeta\in \bC$, $|\zeta|<1$,  
and hence $q_2|_{\IIm(f)}$ is
the identity map.
Now it suffices for us to let $p_4=q_2\circ q_1.$  

This concludes the discussions on embeddings of the simple types.  We now combine the results from the above discussions with those from Proposition 2.  According to Proposition 2a, the representation $\rho$ involved is of form
$$\sum_{i=1}^{r_1}\rho_i+\sum_{i=r_1+1}^{r_1+r_2}\iota_{p_i,q_i}\circ\rho_i.$$  

Suppose $(r_1,r_2)=(1,0)$ or $(0,1)$, that is the representation is irreducible.
Then the projection occurs from composition of projections corresponding to Type 1, 2, 3 embeddings respectively, making use of Proposition 2a.

Consider now the general case.  Note that each of the factors above corresponds to image lying in some Type III classical domain
$III_{s_i}$ for $1\leqslant i\leqslant r_1$ or $III_{p_j+q_j}$ for $r_1+1\leqslant j\leqslant r_1+r_2$.  As explained in \S2 of [Sat],
this corresponds to diagonal blocks of square matrices in $III_g$, where $\sum_{i=1}^{r_1}s_i+\sum_{j=r_1+1}^{r_1+r_2}(p_i+q_i)\leqslant g$.
For simplicity of notation, let us just define $s_j=p_j+q_j$ for $r_1+1\leqslant j\leqslant r_1+r_2$.
Hence we have 
$$i:I_{p,1}\stackrel{i_a}\longrightarrow \left(\begin{array}{cccc}
III_{s_1}&0&&\cdots\\
0&III_{s_2}&&\cdots\\
\vdots&\vdots&\dots&
\end{array}\right)\stackrel{i_b}\hooklongrightarrow I_g.
$$
The projection 
$$p_b:I_g\longrightarrow \left(\begin{array}{cccc}
III_{s_1}&0&&\cdots\\
0&III_{s_2}&&\cdots\\
\vdots&\vdots&\dots&
\end{array}\right)$$
so that $p_b\circ i_b$ is identity is constructed exactly as in those Type $i_1$ embeddings discussed earlier.  For each $III_{s_i}, 1\leqslant i\leqslant r_1+r_2$, there
is a projection $p_{s_i}:III_{s_i}\rightarrow I_{p,1}$ so that $p_{s_i}\circ i_1$ is the identity map according to earlier discussions.  The projection 
$$p_a:\left(\begin{array}{cccc}
III_{s_1}&0&&\cdots\\
0&III_{s_2}&&\cdots\\
\vdots&\vdots&\dots&
\end{array}\right)\longrightarrow I_{p,1}$$
is then defined by identifying $I_{p,1}$ with $i_a(I_{p,1})$ and 
letting 
$$p_a=\frac1{r_1+r_2}p_{s_1}+\cdots + \frac1{r_1+r_2}p_{s_{r_1+r_2}},$$
where addition is given in terms of the coordinate functions of the standard realization, which is the Harish-Chandra coordinates for $I_g$.
Clearly from construction, $p_a\circ i_a$ is the identity.  We may then define $p=p_a\circ p_b$.  It follows from construction that $p\circ i$ is the identity.
 \qed

\bs
\noi{\bf 3.5}  Though not really needed for this article, we mention that the argument of Proposition 3 can 
be  applied to other pairs of Hermitian symmetric spaces of non-compact type, following case-by-case checking as done above for classical domains
using the results of [Sat], see also [I].  In a more inspiring way, Ngaiming Mok [Mok3] has shown us a conceptual proof of such a result for all Hermitian symmetric spaces of non-compact type,
including those containing factors of exceptional types without using classification results.
This was done in terms of the Lie triple system and Harish-Chandra embedding for all pairs of Hermitian symmetric spaces of non-compact type.

\bs

\begin{center}
{\bf \S4. Proofs of the main results}
\end{center}

\bs
\noi{\bf 4.1} Denote by $g_{V,K}$ the Kobayashi (pseudo-)metric of a variety $V$, which is the positive semi-definite Finsler metric defined as
$$\sqrt{g_{V,K}}(x,v)=\inf\{\frac1r|\exists f:\Delta_r\rightarrow V \ {\mbox{holomorphic}}, f(0)=x, f'(0)=v\},$$
where $\Delta_r$ is the disk of radius $r$ in $\bC$ centered at the origin.  It follows from definition that 
the metric on a manifold is the same as on its universal covering from the lifting properties of a map from the unit
disk.  Note for orbifolds, we are considering orbifold maps and orbifold uniformization as explained in {\bf 2.1}.  

From Ahlfors Schwarz Lemma, it follows easily that the Kobayashi metric on a complex ball is precisely
the same as the Poincar\'e metric which is the same as the Bergman metric.  The reader may consult Proposition 3 of [Y2] and the references quoted there
for various forms of Schwarz Lemma.  Since we are considering quotients of bounded domains, it follows from Schwarz Lemma 
that the Kobayashi metric is positive definite in this paper.  Note that the Teichm\"uller space $\cT_g$ can be realized as a bounded
domain in $\bC^{3g-3}$ from Bers Embedding.

Furthermore, it follows immediately from the definition that $g_{V,K}$ has decreasing 
properties in the following sense.  Let $F:M\rightarrow N$ be a holomorphic mapping.  Then
$g_{N,K}(F_*v)\leqslant g_{M,K}(v)$ for all $v\in T_xM$.
Again, the Kobayashi metric may be degenerate in general, but in our case it is always non-degenerate.  This follows from the fact
that the manifolds involved by uniformized by bounded domains in $\bC^n$ for some $n>0$ and the earlier discussions.

We have the following consequence of the decreasing property of the Kobayashi metric.

\begin{lemm} 
Suppose $M=B^N_{\bC}/\Gamma$ is a totally geodesic subvariety of $\cA_g$.  Then  $g_{M,K}=(j_g^{-1})^*g_{T_g^o,K}|_M=g_{\cA_g,K}|_M$.
\end{lemm}

\noi{\bf Proof} The inclusion map $i:M\rightarrow \cA_g$ is a holomorphic embedding.
  Now we have the holomorphic mappings
\begin{equation}
M\rightarrow \cM_g^o \rightarrow \cA_g.
\end{equation}
The first holomorphic map in (12) comes from our assumption that $M\subset T_g^o$ and the fact that $j_g^{-1}|_{T_g^o-j_g(H_g)}$ is a holomorphic map.  Here we used 
the fact that $j_g$ is an injective holomorphic map and is an immersion on $\cM_g-H_g$.  
The second holomorphic mapping in (12)
follows from Torelli mapping.
 The Torelli mapping is holomorphic by 
definition.
Since $M$ is a complex submanifold of  $T_g^o\subset \cA_g$, it follows from definition of the Kobayashi metric in terms of extremal functions that

\begin{equation} g_{M,K}\geqslant (j_g^{-1})^*g_{T_g^o,K}|_M\geqslant g_{\cA_g,K}|_M.
\end{equation}

On the other hand, the Kobayashi metric on a manifold $g_{V,K}$ is the same as its lift $g_{\tV,K}$
to the universal covering $\tV$ of $V$.  Hence in terms of a
 totally geodesic $B^N$ in Siegel $\cS_g$, we need to compare $g_{B^N,K}$ and $g_{\cS_g,K}|_{B^N}$.
  It follows  from Proposition 3 that there is a holomorphic map
 $p:\cS_g\rightarrow \tM$ so that $p\circ i$ is identity.  Hence for $x\in \tM\subset \cS_g$ and $w\in T_x\tM\subset T_x\cS_g$, a
 holomorphic curve  $f:\Delta_r\rightarrow \cS_g \ {\mbox{holomorphic}}$ with $f(0)=x, f'(0)=w$ gives rise to a holomorphic map $p\circ f:\Delta_r\rightarrow \tM$ holomorphic with $p\circ f(0)=x$ and $(p\circ f)'(0)=w$.
 It follows from the decreasing property of the Kobayashi metric
 that $g_{\cS_g,K}|_{B^N}\geqslant  g_{B^N,K}$, which is equivalent to $g_{\cA_g,K}|_{M}\geqslant g_{M,K}$ after descending to $M$ from the universal covering
 as discussed earlier.
 
 Combining the above two paragraphs, we conclude that $g_{\cS_g,K}|_{B^N}=g_{B^N,K}$  and that the two inequalities in (13) can be replaced by equalities.
 
\qed

\bs
\noi{\bf 4.2} Theorem 1 now follows by putting the earlier arguments together.

\ms\noi
{\bf Proof of Theorem 1}  
 Assume for the sake of proof by contradiction that 
there exists $M=B^N_{\bC}/\Gamma$ so that $M$ is a totally geodesic subvariety of $\cA_g$ with $N>1$, and $M\subset T_g^o-j_g(H_g)=j_g(\cM_g^o),$
where $\cM_g^o=\cM_g-H_g$.  

From  Lemma 2, $g_{M,K}=(j_g^{-1})^*g_{T_g^o,K}|_M=g_{\cA_g,K}|_M$.  
In particular,  $(j_g^{-1})^*g_{\cM_g^o,K}|_M=g_{M,K}$.  This however contradicts Proposition 1.  

\qed

\bs
\noi{\bf 4.3}  We remark that the argument in the proof of Theorem 1 can be applied to study the non-existence of locally Hermitian symmetric space in $T_g^o- j_g(H_g)$
as a totally geodesic complex suborbifold for $g\geqslant 2$, a result 
proved earlier in [Ha].      For this purpose, we observe that an analogue of Proposition 1 is true for $\tM$ being
 any Hermitian symmetric space as given in [A].   Together with the remarks given in {\bf 3.5} and the results of [Mok3], the other parts of the proof can be applied.
 

\bs
 \noi{\bf 4.4} {\bf Proof of Theorem 2}

From the results of [Ha], we know that any symmetric variety $M$ in 
$T_g^o- j_g(H_g)\subset \cA_g$ has to be of real rank $1$ as a locally symmetric space.  
Alternatively, to make the article more self-contained, this follows from {\bf 3.5, 4.3} and the proof of Theorem 1.
Since $M$ is Hermitian symmetric, we know that $M$ has to be a complex ball quotient.
Theorem 1 implies that the complex dimension of $M$ is $1$ and hence $M$ is a hyperbolic Riemann surface.   From the discussions above, such a Riemann surface $M$
has to be a Shimura-Teichm\"uller in the terminology of [M\"o], since the Kobayashi metric, which is well-known to be the same as the Teichm\"uller metric on Teichm\"uller spaces, is the same as the natural hyperbolic metric
on $M$ as it is a totally geodesic curve in $T_g^o- j_g(H_g)$.
In such a case, M\"oller proved in [M\"o]
that such a Riemann surface does not exist for genus $g\neq 3,4,5$, and the only examples for $g=3,4$ are given in the statement of Theorem 2.   Very recently, it was
proved by Aulicino and Norton in [AN] that there is no example in $g=5$.  Theorem 2 follows.
\qed

\bs
Theorem 2 gives a necessary and sufficient condition for the existence of a locally Hermitian symmetric space in $T_g^o- j_g(H_g)$.

\bs
\noindent{\bf References} 

\bs

\ms
\noi [A]
Antonakoudis, S. M.,
Teichm\"uller spaces and bounded symmetric domains do not mix isometrically. 
Geom. Funct. Anal. 27 (2017), no. 3, 453-465.

\ms
\noi[AN] Aulicino, D., Norton, C., Shimura-Teichm\"uller curves in genus $5$, arXiv:1903.01625v2.  

\ms
\noi [CLTZ] Chen, K., Lu, X., Tan, S-L., Zuo, K., On Higgs bundles over Shimura varieties of ball quotient type. Asian J. Math. 22 (2018), no. 2, 269-284.

\ms
\noi
[C] Coleman, R. F., Torsion points on curves, in Galois representations and arithmetic algebraic geometry (Kyoto, 1985/Tokyo, 1986), Advanced Studies in Pure Mathematics, vol. 129 (North-Holland, Amsterdam, 1987), 235-247.

\ms
\noi [DM] Daskalopoulos, G.,  Mese, C., Rigidity of Teichm\"uller space, arXiv:1502.03367.

\ms\noi
[FM] Farb B., and Masur, H., Superrigidity and mapping class groups, Topology 37 (1998), 
1169-1176.

\ms\noi
[dJN] de Jong, A., Noot, R., Jacobians with complex multiplication, in: Arithmetic Algebraic Geometry, Texel, 1989, in: Progr.
Math., vol. 89, Birkh\"auser, 1991, pp. 177-192.

\ms\noi
[dJZ] de Jong, A., Zhang, S., Generic abelian varieties with real multiplication are not Jacobians. Diophantine geometry, 165-172, CRM Series, 4, Ed. Norm., Pisa, 2007.

\ms
\noi [Ha] Hain, R., Locally symmetric families of curves and Jacobians. In: Moduli of curves and abelian varieties; pp. 91-108. (C. Faber, E. Looijenga, eds.) Aspects of Math. E33, Vieweg, Braunschweig, 1999.

\ms
\noi [HC] Harish-Chandra, Representations of semisimple Lie groups, CI, Amer. J. Math. 78 (1956), 564-628.

\ms
\noi [He] Helgason, S., Differential geometry, Lie groups, and symmetric spaces, Pure Appl.
Math., 80, Academic Press, 1978.


\ms
\noi [I] Ihara, S., Holomorphic imbeddings of symmetric domains. J. Math. Soc. Japan 19(1967), 261-302.

\ms
\noi [KS] Korevaar, N., Schoen, R., Global existence theorem for harmonic maps to non-locally compact spaces. Comm. Anal. Geom. 5 (1997) 333-387.

\ms
\noi
[Mok1] Mok, N., 
Metric Rigidity Theorems On Hermitian Locally Symmetric Manifolds, World Scientific, 1989.

\ms
\noi
[Mok2] Mok, N., Projective algebraicity of minimal compactifications of complex-hyperbolic space forms of finite volume. Perspectives in analysis, geometry, and topology, 331?354, Progr. Math., 296, Birkh\"auser/Springer, New York, 2012.

\ms
\noi
[Mok3] Mok, N., private communication.

\ms
\noi [MSY] Mok, N., Siu, Y.-T., Yeung, S.-K., Geometric superrigidity,  Invent. Math., 113(1993), 57-84.


\ms
\noi [M\"o] M\"oller, M., Shimura and Teichm\"uller curves. J. Mod. Dyn. 5 (2011), 1-32.

\ms
\noi [Moo]  Moonen, B., Special subvarieties arising from families of cyclic covers of the projective line, Doc. Math. 15 (2010) 793-819.

\ms
\noi [O] Oort, F., Canonical liftings and dense sets of CM-points, in Arithmetic geometry (Cortona, 1994), Symposia Mathematica, vol. XXXVII (Cambridge University Press, Cambridge, 1997), 228-234.

\ms
\noi [Sag]  Sagle, A. A.,  A note on triple systems and totally geodesic submanifolds in a homogeneous space. Nagoya Math. J. 32 (1968), 5-20.


\ms
\noi [Sat] Satake, I., Holomorphic imbeddings of symmetric domains into a Siegel space. Amer. J. Math. 87(1965), 425-461.


\ms
\noi [Si1] Siu, Y.-T.,  The complex-analyticity of harmonic maps and the strong
rigidity of compact K\"ahler manifolds, Ann. of Math. 112 (1980), 73-111.

\ms
\noi [Si2] Siu, Y.-T., Strong rigidity of compact quotients of exceptional bounded symmetric domains. Duke Math. J. 48 (1981), 857-871.

\ms
\noi [SY] Siu, Y.-T., Yau, S.-T., Compactification of negatively curved complete K\"ahler manifolds of finite volume. Seminar on Differential Geometry, pp. 363-380, 
Ann. of Math. Stud., 102, Princeton Univ. Press, Princeton, N.J., 1982.




\ms
\noi [Y1] Yeung, S.-K., Uniformization of 1/4-pinched negatively curved manifolds with special holonomy. Int. Math. Res. Not. 1995, 365-375.

\ms
\noi [Y2] Yeung, S.-K.,  Quasi-isometry of metrics on Teichm\"uller spaces. Int. Math. Res. Not. 2005, 239-255. 
\end{document}